\documentclass{amsart}
\usepackage[]{latexsym,amssymb,amsmath,amsfonts, amsthm, verbatim}
\usepackage[all,cmtip,ps]{xy}
\newtheorem{TEO}{Theorem}[section]
\newtheorem{PROP}[TEO]{Proposition}
\newtheorem{COR}[TEO]{Corollary}
\newtheorem{LEMMA}[TEO]{Lemma}
\theoremstyle{definition}
\newtheorem{DEF}[TEO]{Definition}
\newtheorem{EX}[TEO]{Example}
\newtheorem{REM}[TEO]{Remark}

\let\End\undefined

\def\Hom{\mathop{{\rm Hom}}}
\def\End{\mathop{{\rm End}}}
\def\R{\mathbf{{{R}}}}
\def\pd{\mathop{{\rm pd}}}

\def\Ext#1#2#3#4{\mathop{{\mathrm {Ext}}^{#1}_{#2}(#3,#4)}}

\renewcommand*{\Phi}{\varPhi}
\renewcommand*{\Psi}{\varPsi}

\newcommand*{\lmod}{\textrm{\textup{-Mod}}}

\def\seb#1#2#3#4#5{0\rightarrow#1\stackrel{#4}{\rightarrow}#2
\stackrel{#5}{\rightarrow}#3\rightarrow 0}

\def\fdim{\mathop{{\rm {\mathcal F}dim}}}
\def\chdim{\mathop{{\rm ch.dim}}}
\def\glfdim{\mathop{{\rm gl{\mathcal F}dim}}}
\def\gldim#1{\mathop{{\rm gl{\mathcal #1}dim}}}
\newcommand{\A}{\mathcal{A}}
\newcommand{\B}{\mathcal{B}}
\newcommand*{\D}{\mathcal{D}}
\newcommand*{\Hm}{H}

\begin{document}
\title{On classes defining a homological dimension}
\author[F. Mantese]{Francesca Mantese}
\author[A. Tonolo]{Alberto Tonolo}
\address[F. Mantese]{Dipartimento di Informatica, Universit\`a degli Studi di Verona, strada Le Grazie  15, I-37134 Verona - Italy}
\email{mantese@sci.univr.it}
\address[A. Tonolo]{ Dip. Matematica Pura ed Applicata, Universit\`a degli studi di Padova, via Trieste 63, I-35121 Padova Italy}
\email{tonolo@math.unipd.it}
\thanks{Research of the second author supported by grant CDPA048343 of Padova University}
\date{}
\begin{abstract}
A class $\mathcal F$ of objects of an abelian category $\mathcal A$ is said to define a \emph{homological dimension} if for any object in $\mathcal A$ the length of any $\mathcal F$-resolution is uniquely determined.  In the present paper we investigate classes satisfying this property. 

\end{abstract}

\maketitle

\section*{Introduction}

In general the class of the objects of a given abelian category $\mathcal A$ is too complex to admit any satisfactory classification. Starting from a known subclass $\mathcal F$ of $\mathcal A$, one may try to approximate arbitrary objects by the objects in $\mathcal F$. This approach has successfully been followed  over the past few decades for categories of modules through the theory of precovers and preenvelopes, or left and right approximations (see \cite{EJ} or \cite{GT} for a detailed list of references).

Another point of view could be to measure the ``distance'' of any object in $\mathcal A$ from the class $\mathcal F$, introducing a notion of \emph{dimension} with respect to the class $\mathcal F$, computed by means of $\mathcal F$-resolutions.  In this framework, the notions of projective dimension, weak dimension, Gorenstein dimension of modules have been deeply studied.

Our aim is to define a good concept of dimension with respect to a wide family of classes of objects. We say that a class $\mathcal F$ of objects of an abelian category $\mathcal A$ defines a \emph{homological dimension} if for any object in $\mathcal A$, the length of any $\mathcal F$-resolution is uniquely determined (see Definition~\ref{definition:homological}). In such a way to each object in $\mathcal A$ one can associate an $\mathcal F$-invariant number which represents locally the relevance of $\mathcal F$.

In the first section we study several properties of classes defining a homological dimension; in particular we discuss their closure properties and the connection with precover classes and cotorsion pairs. 
In the second section, using tools from derived categories, we generalize the Auslander notion of Gorenstein dimension to arbitrary abelian categories. We consider a homological dimension associated to an adjoint pair $(\Phi,\Psi)$ of contravariant functors, obtaining again the classical Gorenstein dimension on $R$-modules in case $\Phi=\Psi=\Hom(-,R)$ for a commutative noetherian ring $R$.

\section{Homological dimension}



\begin{DEF}[conf. \cite{AB}] Let $\mathcal F$ be a class of objects in an abelian category $\mathcal A$. We say that an object $M$ in $\mathcal A$ has \emph{left $\mathcal F$-dimension $\leq \alpha$}, $\alpha\in\mathbb N\cup \{\infty\}$, if there exists a long exact sequence
\[ ...\to F_i\to F_{i-1}\to ...\to F_1\to F_0\to M\to 0\]
with $F_i\in \mathcal F\cup\{0\}$, and $F_i=0$ for $i>\alpha$. We denote by $\mathcal F_{\alpha}$, the class of objects $M$ of left $\mathcal F$-dimension $\leq \alpha$ (shortly $\fdim M\leq\alpha$), and by $\mathcal F_{<\infty}$ the class of objects of finite left $\mathcal F$-dimension.
\end{DEF}

In general there exist objects which have not a left $\mathcal F$-dimension: in particular all objects which are not quotients of objects in $\mathcal F$.  
We denote by $\overline{\mathcal F}$ the class of all objects in $\mathcal A$ which are homomorphic image of objects in $\mathcal F$.

\begin{REM}
If $\mathcal A$ has enough projectives and $\mathcal F$ is closed under direct summands, then $\overline{\mathcal F}=\mathcal A$ if and only if $\mathcal F$ contains all projective objects. 

In particular, if $\mathcal A=R\lmod$, denoted by $\mathcal P$ and $\mathcal Fl$ the classes of projective and flat modules respectively, then $\overline{\mathcal P}=R\lmod$ and $\overline{\mathcal Fl}=R\lmod$, and left $\mathcal P$- and left $\mathcal Fl$- dimensions are the usual projective and flat (or weak) dimensions of a module.
\end{REM}

\begin{DEF}
We say that $\mathcal A$ has \emph{global left $\mathcal F$-dimension} $\leq\alpha$ (resp. $<\infty$), $\alpha\in\mathbb N\cup \{\infty\}$, if for each object $M$ in $\mathcal A$ we have $\fdim M\leq\alpha$ (resp. $<\infty$).
\end{DEF}
Clearly $\mathcal A$ has global left $\mathcal F$-dimension $\leq\infty$ if and only if $\mathcal A=\overline{\mathcal F}$.

In any abelian category $\mathcal A$ it is possible (see \cite[Ch. VII]{Mi}) to define, for any pair of object $A$, $B\in\mathcal A$, the family $\Ext i{\mathcal A}AB$ of equivalence classes of exact sequences of length $i$ with left end $B$ and right end $A$, with respect to the Yoneda equivalence relation. The family $\Ext i{\mathcal A}AB$ in general is not a set (see \cite[Ch. VI]{F}); nevertheless it can be equipped with an additive structure and become a \emph{big abelian group}. The big abelian groups are defined in the same way as ordinary abelian groups, except than the underlying class need not be a set. Quoting \cite{Mi}, ``[...] we are prevented from talking about  the category of big abelian groups because the class of morphisms between a given pair of big groups need not be a set.  Nevertheless this will not keep us from talking about kernels, cokernels, images, exact sequences, etc., for big abelian groups.'' If $\mathcal A$ has enough injectives or projectives, then $\Ext i{\mathcal A}AB$ is an abelian group for each $A$, $B\in\mathcal A$.

Given a class of objects $\mathcal G$, we denote by
\[\mathcal G^{\perp_m}=\{M\in\mathcal A:\Ext i{\mathcal A}GM=0,\ \forall 1\leq i\leq m,\ G\in\mathcal G\};\]
the intersection $\bigcap_{m\geq 1}\mathcal G^{\perp_m}$ will be denoted by $\mathcal G^{\perp_{\infty}}$.
Dually, we denote by
\[{}^{\perp_m}\mathcal G=\{M\in\mathcal A:\Ext i{\mathcal A}MG=0,\ \forall 1\leq i\leq m,\ G\in\mathcal G\};\]
the intersection $\bigcap_{m\geq 1} {}^{\perp_m}\mathcal G$ will be denoted by $^{\perp_{\infty}}\mathcal G$.

\begin{DEF}[\text{\cite[Ch. VI.6]{Mi}}] Let $A$ be an object of an abelian category $\mathcal A$. The \emph{cohomological dimension $\chdim A$ of $A$} is the least integer $n$ such that the one variable functor $\Ext n{\mathcal A}-A$ is not zero.
\end{DEF}

If $\mathcal A$ has enough injective objects (e.g., if $\mathcal A$ is a Grothendieck category) the cohomological dimension of an object coincides with its injective dimension. 

\begin{PROP}\label{prop:gldim}
Assume that $\mathcal A$ has enough projectives.
\begin{enumerate}
\item If $\glfdim \mathcal A\leq n$, $n\in\mathbb N$, then $\chdim Y\leq n$ for each $Y\in\mathcal F^{\perp_{n+1}}$.  
\item If $\mathcal F= {^{\perp_m}{\mathcal G}}$ for a class $\mathcal G$ of modules of cohomological dimension less or equal than  $n\in\mathbb N$, then $\glfdim \mathcal A\leq n$.
\end{enumerate}
\end{PROP}
\begin{proof}
Let $M$ be an arbitrary object in $\mathcal A$.

1)  Since $\glfdim \mathcal A\leq n$, there exists an exact sequence
\[0\to F_n\to F_{n-1}\to \dots\to F_1\to F_0\to M\to 0.\]
Applying the contravariant functor $\Hom(-,Y)$, since $Y\in\mathcal F^{\perp_{n+1}}$, by dimension shift we get
$\Ext{n+1}{\mathcal A}MY\cong \Ext 1{\mathcal A}{F_n}Y=0$. Since $\Ext{n+1}{\mathcal A}MY=0$ for each object $M$ in $\mathcal A$, and  the latter has enough projectives, then $\Ext{n+i}{\mathcal A}MY=0 $ for each $i\geq 1$, i.e. $\chdim Y\leq n$.

2) Consider an exact sequence
\[0\to K_n\to P_{n-1}\to \dots\to P_1\to P_0\to M\to 0\] 
with $P_i$ projective for $i=0,\dots n-1$. Since $P_i\in\mathcal F$ it is enough to prove that $K_n$ belongs to $\mathcal F$. So let $G\in\mathcal G$; then $\Ext i{\mathcal A} {K_n} G\cong\Ext {n+i}{\mathcal A} M G=0$ for $1\leq i$. Therefore $K_n\in {^{\perp_{\infty}}{\mathcal G}} \subseteq{^{\perp_m}{\mathcal G}}=\mathcal F$.
\end{proof}

In order to introduce a good measure of the  distance between an object of $\A$ and a given class $\mathcal F$, the length of a $\mathcal F$-resolution has to be uniquely determined.

\begin{DEF}\label{definition:homological} We say that the left $\mathcal F$-dimension associated to a class $\mathcal F$ is \emph{ homological} (or that the class $\mathcal F$ defines a \emph{homological dimension}) if
\begin{enumerate}
\item for any short exact sequence $0\to K\to F\to M\to 0$ with $F\in\mathcal F$ and $M\in\mathcal F_{\infty}$, the object $K$ belongs to $\mathcal F_{\infty}$;
\item
for any exact sequence
\[0\to K_n\to F_{n-1}\to \dots\to F_1\to F_0\to X\to 0\]
with $F_i\in\mathcal F$, $i=0,1,...,n-1$, and $X\in\mathcal F_n$,  the object $K_n$ belongs to $\mathcal F$.
\end{enumerate}
\end{DEF}
Clearly  if $\mathcal A=\overline{\mathcal F}$ we have $\mathcal A=\mathcal F_{\infty}$ , and the first condition is empty.

\begin{EX}
If $\mathcal A=R\lmod$, the classes  $\mathcal P$ and $\mathcal Fl$ define a homological dimension.  The class of free modules defines a homological dimension if and only if it coincides with the class of projective modules (see Proposition~\ref{prop:chiusaaddendi}), e.g. if $R$ is local. 

If $\mathcal A$ is the category of coherent sheaves on a noetherian scheme $X$, the classes of the locally free sheaves $\mathcal {LF}$ and of the invertible sheaves $\mathcal {I}$ both define a homological dimension (see \cite[Chp. 2 \S 5, Chp.3 \S 6]{H}). If $X$ is quasi-projective over $\text{Spec} R$, where $R$ is a noetherian commutative ring, then $\overline{\mathcal {LF}}=\mathcal A$.
\end{EX}

Note that the notion of homological dimension can be easily dualized obtaining a notion of \emph{homological codimension}; for instance, if $\mathcal A=R\lmod$, the class  $\mathcal I$ of injective modules defines a homological codimension. Most of the results we obtain in this paper could be reformulated for this dual concept.

In the sequel we study closure properties of classes defining a homological dimension.

Let $\mathcal F$ be a class of modules and $\seb AFC{}{}$ be an exact sequence with $F\in\mathcal F$. Thus, for any $i\geq 1$ in $\mathbb N$, if $A\in\mathcal F_{i-1}$ then $C\in\mathcal F_i$. 
\begin{LEMMA}\label{lemma:kernel}
Let $\mathcal F$ be a class of objects in $\mathcal A$ and $\seb AFC{}{}$ be an exact sequence with $F\in\mathcal F$. If $\mathcal F$ defines a homological dimension and $C\in\mathcal F_i$, then $A\in\mathcal F_{i-1}$. In particular $\mathcal F$ is closed under kernels of epimorphisms.
\end{LEMMA}
\begin{proof}
By the definition of homological dimension, $A$ belongs to $\mathcal F_{\infty}$. Therefore consider an exact sequence \[0\to K_{i-1}\to F_{i-2}\to \dots\to F_1\to F_0\to A\to 0\] with $F_j\in\mathcal F$. Since 
 \[0\to K_{i-1}\to F_{i-2}\to \dots\to F_1\to F_0\to F \to C\to 0\] is an $\mathcal F$-resolution for $C$ and $\fdim C\leq i$, we get that $K_{i-1}\in\mathcal F$.
  \end{proof}

\begin{PROP}\label{prop:chiusaaddendi}
Let $\mathcal F$ be a class of objects defining a homological dimension. If $\mathcal F$ is closed under countable direct sums, then 
$\mathcal F$ is closed under direct summands.
\end{PROP}
\begin{proof}
Let $L\oplus M=F\in\mathcal F$; consider the short exact sequence
\[0\to L\to L\oplus( M\oplus L)^{(\omega)}\to (M\oplus L)^{(\omega)}\to 0;\]
since both $(M\oplus L)^{(\omega)}$ and $L\oplus( M\oplus L)^{(\omega)}\cong (L\oplus M)^{(\omega)}$ belong to $\mathcal F$, also $L$ belongs to $\mathcal F$.
\end{proof}

 In the next theorem we compare the $\mathcal F$-dimension of objects in a short exact sequence.

\begin{TEO}\label{prop:ses} Assume $\mathcal F$ defines a homological dimension and it is closed under finite direct sums.
Let $\seb ABC{}{}$ be a short exact sequence. Then for each $i\in\mathbb N$ we have that
\begin{enumerate}
\item[$(1_i)$] if $B$ and $C$ belong to $\mathcal F_i$ then $A$ belongs to $\mathcal F_i$;
\item[$(2_i)$]  if $A$ and $B$ belong to $\mathcal F_i$ then $C$ belongs to $\mathcal F_{i+1}$.  
\end{enumerate}
If $\overline{\mathcal F}$ is closed under extensions, then
\begin{enumerate}
\item[$(3_i)$] if $A$ and $C$ belong to $\mathcal F_{i+1}$, then $B$ belongs to $\mathcal F_{i+1}$;
\item[$(4_i)$] if $B\in\mathcal F_i$ and $C\in\mathcal F_{i+1}$, then $A$ belongs to $\mathcal F_{i}$.
\end{enumerate}
\end{TEO}
\begin{proof}
(1) - (2): If $i=0$, $2_0$ is clearly true by definition and $1_0$ follows by $\fdim C=0\leq 1$ and the fact that $\mathcal F$ defines a homological dimension. Assume $1_{i-1}$ and $2_{i-1}$ true for ${i-1}\geq 0$.
Let us consider the pullback diagram
$$(*)\quad\xymatrix{& 0 & 0\\
0 \ar[r] & A\ar[r]\ar[u] & B\ar[r]\ar[u] & C\ar[r] & 0\\
0 \ar[r] & P_B\ar[r]\ar[u] & F_B \ar[r]\ar[u] & C\ar[r]\ar@{=}[u] & 0\\  
& K_B\ar[u]\ar@{=}[r] & K_B \ar[u]\\
& 0\ar[u] & 0\ar[u]}$$
with $F_B$ in $\mathcal F$.

$1_i$: 
by Lemma~\ref{lemma:kernel} both $K_B$ and $P_B$ in diagram $(*)$ belong to $\mathcal {F}_{i-1}$, and so by induction $A\in\mathcal F_i$.

$2_i$: Let now $A$ and $B$ be in $\mathcal F_i$; there exist $F_B\in\mathcal F$ and an epimorphism $\pi:F_B\to B$. Consider the following pullback diagram
$$\xymatrix{&& 0 & 0\\
0 \ar[r] & A\ar[r]& B\ar[r]^g\ar[u] & C\ar[r] \ar[u] & 0\\
0 \ar[r] & A\ar[r]\ar@{=}[u] & P_C\ar[r]^p\ar[u] & F_B\ar[r]\ar[u]_{g\circ \pi} & 0\\  
&& K_C\ar[u]\ar@{=}[r] & K_C\ar[u]\\
&& 0\ar[u] & 0\ar[u]}$$
Since $P_C$ is a pullback, there exists $j:F_B\to P_C$ such that $p\circ j=1_{F_B}$.
Then the middle exact sequence splits, and therefore $P_C=A\oplus F_B$; since $\mathcal F$ is closed under finite direct sums, $P_C$ belongs to $\mathcal F_i$. Therefore by $1_i$ we have $K_C\in\mathcal F_i$ and hence $C$ belongs to $\mathcal F_{i+1}$.

(3) - (4): If $i=0$, $4_0$ follows by the definition of homological dimension. Since $\overline{\mathcal F}$ is closed under extensions, if $A$ and $C$ are in $\overline{\mathcal F}$, also $B$ belongs to $\overline{\mathcal F}$. Then, if $A$ and $C$ belong to $\mathcal F_1$, we can consider the pullback diagram $(*)$ with $F_B$ in $\mathcal F$. Since $C$ belongs to $\mathcal F_1$, then $P_B$ belongs to $\mathcal F$; since $A$  belongs to $\mathcal F_1$, then also $K_B$ belongs to $\mathcal F$, and therefore $B$ belongs to $\mathcal F_1$.
Assume $3_{i-1}$ and $4_{i-1}$ true for ${i-1}\geq 0$. 

$4_i$: Let us consider the pullback diagram
$$\xymatrix{& & 0 & 0\\
0 \ar[r] & A\ar[r] & B\ar[r]\ar[u] & C\ar[r]\ar[u] & 0\\
0 \ar[r] & A\ar[r]\ar@{=}[u] & P_C\ar[r]\ar[u] & F_C\ar[r]\ar[u] & 0\\  
& & K_C\ar[u]\ar@{=}[r] & K_C\ar[u]\\
& & 0\ar[u] & 0\ar[u]}$$
with $F_C\in\mathcal F$; then $K_C$ belongs to $\mathcal F_i$. Since $B$ belongs to $\mathcal F_i$, by $3_{i-1}$ we have that $P_C\in \mathcal F_i$, and hence, by $1_i$, $A$ belongs to $\mathcal F_i$.

$3_i$: Since $\overline{\mathcal F}$ is closed under extensions, we can consider the pullback diagram $(*)$ with $F_B$ in $\mathcal F$. By Lemma~\ref{lemma:kernel}, $P_B$ belongs to $\mathcal F_i$; then $K_B\in\mathcal F_i$ by $4_i$, and hence $B$ belongs to $\mathcal F_{i+1}$.
\end{proof}

\begin{REM}
It follows that if $\mathcal F$ is closed under finite direct sums and $\overline{\mathcal F}$ is closed under extensions, then
\begin{itemize}
\item the class $\mathcal F_{<\infty}$ is closed under extensions, kernels of epimorphisms and cokernels of monomorphisms;
\item the classes $\mathcal F_i$, $i\geq 0$, are closed under  kernels of epimorphisms; if $i\geq 1$, they are closed also under extensions.
\end{itemize}
\end{REM}

\begin{PROP} Assume $\mathcal F$ defines a homological dimension, it is closed under finite direct sums, and $\overline{\mathcal F}=\mathcal A$.
Then also $\mathcal F_i$ and $\mathcal F_{<\infty}$ define a homological dimension for any $i\geq 1$. 
\end{PROP}
\begin{proof}
Since $\overline{\mathcal F}=\mathcal A$, also $\overline{\mathcal F_i}=\mathcal A=\overline{\mathcal F_{<\infty}}$. Therefore condition 1 in Definition~\ref{definition:homological} is empty in both the cases.
Let $M$ be an object admitting an $\mathcal F_i$-resolution
$$ 0\to F_{i,n}\to F_{i, n-1}\to \dots\to F_{i,0}\to M\to 0 .$$
Consider an exact sequence  $ 0\to K \to F'_{i,n-1}\to \dots\to F'_{i,0}\to M\to 0 $ with $F'_{i,j}\in\mathcal F_i$. From the first sequence, applying recursively Theorem~\ref{prop:ses}, 2), we get that $M\in\mathcal F_{n+i}$. Applying recursively Theorem~\ref{prop:ses}, 4) to the second exact sequence we obtain that $K\in\mathcal F_i$. 
Since each finite $\mathcal F_{<\infty}$ resolution is actually an $\mathcal F_m$ resolution for a suitable $m\in\mathbb N$, we conclude that also $\mathcal F_{<\infty}$ defines a homological dimension.
\end{proof}

In case the abelian category $\mathcal A$ has enough projectives, a relevant family of classes defining a homological dimension is given by the left orthogonal of any class.

\begin{PROP}\label{prop:perp}
Assume $\mathcal A$ has enough projectives, and let $\mathcal G$ be a class of objects in $\mathcal A$.
 Then $\mathcal F={}^{\perp_m}{\mathcal G}$, $1\leq m\in\mathbb N$, defines a homological dimension if and only if \[\mathcal F={}^{\perp_{\infty}}{\mathcal G}.\]
 In such a case $\mathcal A=\overline{\mathcal F}$.
\end{PROP}
\begin{proof}
Assume $\mathcal F={}^{\perp_m}{\mathcal G}$ defines a homological dimension. Let us prove that $\mathcal F={}^{\perp_{m+1}}{\mathcal G}$; then we conclude inductively. Consider an arbitrary object $F\in  \mathcal F$. Consider a short exact sequence
\[0\to K\to P\to F\to 0\]
with $P$ projective; since $P$ belongs to $\mathcal F$, by Lemma~\ref{lemma:kernel} we have that also $K\in\mathcal F$. Therefore for each $G\in\mathcal G$ we have
\[\Ext{m+1}{\mathcal A}FG\cong \Ext{m}{\mathcal A}KG=0,\]
because $K\in\mathcal F$.

Conversely, let us prove that $\mathcal F={}^{\perp_{\infty}}{\mathcal G}$ defines a homological dimension.
Clearly, containing $\mathcal F$ the projectives, each object has left $\mathcal F$-dimension $\leq\infty$.
Let $M$ be an object with $\fdim M\leq n$, $n\in\mathbb N$. Then there exists an exact sequence
\[0\to F'_n\to F'_{n-1}\to \dots\to F'_1\to F'_0\to M\to 0\]
with $F'_i\in\mathcal F$ for $i=0, \dots, n$. Let us consider an exact sequence 
\[0\to K_n\to F_{n-1}\to \dots\to F_1\to F_0\to M\to 0\]
with $F_i\in \mathcal F$ for $i=0,\dots n-1$. Let us show that $K_n\in\mathcal F$. In fact, let $X\in\mathcal G$. Then $\Ext i{\mathcal A}{K_n} X\cong \Ext{n+i}{\mathcal A}M X\cong \Ext i{\mathcal A}{F'_n} X=0$ for each $i\geq 1$.\end{proof}

\begin{EX}\label{EX:ex}
\begin{enumerate}
\item
Since $\mathbb Z$ has global dimension $1$, the class $\mathcal W={^{\perp_1}{\mathbb Z}}={^{\perp_{\infty}}{\mathbb Z}}$ of Whitehead abelian groups defines a homological dimension. By Proposition~\ref{prop:gldim}, \emph{(ii)} we have $\gldim W \mathbb Z\leq 1$.
\item
Any torsion free class in a category of modules defines a homological dimension, since it is closed under submodules. In general it is not the left orthogonal of any class. Consider for example the class $\mathcal R$ of reduced abelian groups; since $\mathcal R^{\perp_{\infty}}$ is the class of divisible groups, ${}^{\perp_{\infty}}(\mathcal R^{\perp_{\infty}})$ is the whole class of abelian groups. Therefore $\mathcal R$ cannot be the left orthogonal of a class, otherwise ${}^{\perp_{\infty}}(\mathcal R^{\perp_{\infty}})$ would be equal to $\mathcal R$.
\end{enumerate}
\end{EX}

In the following results we are interested in giving necessary or sufficient conditions for a class defining a homological dimension to be a left orthogonal.  

\begin{LEMMA}\label{lemma:perps}
Assume $\mathcal A$ has enough projectives. If $\mathcal F$ defines a homological dimension and it contains the projectives, then $\mathcal F^{\perp_1}=\mathcal F^{\perp_{\infty}}$.
\end{LEMMA}
\begin{proof}
Let $M$ be an object in $\mathcal F^{\perp_1}$ and $F\in\mathcal F$.
Consider a short exact sequence $\seb {F'}PF{}{}$ with $P$ projective; since $\mathcal F$ defines a homological dimension also $F'$ belongs to $\mathcal F$. Applying $\Hom_{\mathcal A}(-,M)$ we get $\Ext{i+1}{\mathcal A}FM\cong\Ext i{\mathcal A}{F'}M$; then $\Ext{2}{\mathcal A}FM=0$ and we conclude by induction.
\end{proof}

\begin{TEO} Assume $\mathcal A$ has enough projectives, and 
let $\mathcal F$ be a special precover class. Then $\mathcal F$ defines a homological dimension if and only if $\mathcal F={^{\perp_{\infty}}({{\mathcal F}^{\perp_{\infty}}})}$.  
\end{TEO}
\begin{proof} 
If $\mathcal F={^{\perp_{\infty}}({{\mathcal F}^{\perp_{\infty}}})}$, by Proposition~\ref{prop:perp} we get that $\mathcal F$ defines a homological dimension. 

Conversely, suppose that $\mathcal F$ defines a homological dimension. 
Let us prove that $\mathcal F={}^{\perp_1}(\mathcal F^{{\perp}_{\infty}})$.  Of course $\mathcal F\subseteq {}^{\perp_1}(\mathcal F^{{\perp}_{\infty}})$. Let now $M\in{}^{\perp_1}(\mathcal F^{{\perp}_{\infty}})$; consider a special $\mathcal F$-precover $0\to K\to F\to M\to 0$. Since by the previous lemma  $K\in \mathcal F^{\perp_1}=\mathcal F^{{\perp}_{\infty}}$, we get $\Ext 1RMK=0$. Since the special precover classes are closed under direct summands \cite[Section 2.1]{GT}, then $M\leq^{\oplus} F$ belongs to $\mathcal F$. Again by Proposition~\ref{prop:perp}  we conclude that $\mathcal F={^{\perp_{\infty}}({{\mathcal F}^{\perp_{\infty}}})}$.
\end{proof}

Most of the examples of classes defining a homological dimension give special precovers. Nevertheless observe that this is not always the case: Eklof and Shela in \cite{ES} proved that, consistently with ZFC,  the class of Whitehead abelian groups, which defines a homological dimension (see Example~\ref{EX:ex}), does not provide precovers. In particular they proved that $\mathbb Q$, which has $\mathcal W$-dimension 1, does not admit $\mathcal W$-precover.

\begin{REM}
If $\mathcal F$ is a special precover class and it defines a homological dimension, then for each module $M$ it is possible to get an $\mathcal F$-resolution
\[\dots\to F_i\to\dots\to F_1\to F_0\to M\to 0\]
such that, denoted by $\Omega_{\mathcal F}^i(M)$ the $i$-th $\mathcal F$ syzygy of $M$, the induced map $F_j\to \Omega_{\mathcal F}^{j-1}(M)$ is a special $\mathcal F$-precover of $\Omega_{\mathcal F}^{j-1}(M)$. Therefore, in such a case our definition of $\mathcal F$-dimension coincides with the definition given by Enochs and Jenda (see  \cite[Definition~8.4.1]{EJ}).
\end{REM}

Other significative  classes defining a homological dimension are those studied by Auslander-Buchweitz in \cite{AB}. In that paper they introduced the notion of \emph{$Ext$-injective cogenerator} for an additively closed exact subcategory $\mathcal F$ of $\mathcal A$: an additively closed subcategory  $\omega\subseteq \mathcal F$ is an Ext-injective cogenerator for $\mathcal F$ if $\omega\subseteq {\mathcal F}^{\perp_{\infty}}$ and for any $F\in \mathcal  F$ there exists an exact sequence $0\to F\to X \to F'\to 0$ where $F'\in\mathcal F$ and $X\in\omega$. 

\begin{PROP}\cite[Propositions~2.1, 3.3]{AB} 
Let  $\mathcal F$ be an  additively closed exact subcategory  of $\mathcal A$ closed under kernels of epimorphisms. If $\mathcal F$ admits an $Ext$-injective cogenerator $\omega$, then $\mathcal F$ defines a homological dimension. Moreover, if any object has finite $\mathcal F$-dimension, then $\mathcal F={^{\perp_{\infty}}\mathcal G}$, where $\mathcal G$ is the class of   objects in $\mathcal A$ of finite $\omega$-dimension.
\end{PROP}

We conclude  this section remarking the connection between classes defining a homological dimension  and  cotorsion pairs in categories of modules. So 
 we assume $\mathcal A=R\lmod$, the category of left $R$-modules over a ring $R$.

\begin{DEF}
Let $\mathfrak A$ and $\mathfrak B$ be two classes of modules. The pair $(\mathfrak A,\mathfrak B)$ is called a \emph{cotorsion pair}
if $\mathfrak A={}^{\perp_1}\mathfrak B$ and $\mathfrak A^{\perp_1}=\mathfrak B$. The pair 
$(\mathfrak A,\mathfrak B)$ is called an \emph{hereditary cotorsion pair} if $\mathfrak A={}^{\perp_{\infty}}\mathfrak B$ or equivalently $\mathfrak A^{\perp_{\infty}}=\mathfrak B$.
\end{DEF}

We  stress that,  by Proposition~\ref{prop:perp}, the hereditary cotorsion pairs are exactly the cotorsion pairs $(\mathfrak A,\mathfrak B)$ such that $\mathfrak A$ defines a homological dimension. 
%

\begin{EX}
Let $R$ be a commutative domain. A module $M$ is \emph{Matlis cotorsion} provided that $\Ext 1RQM=0$, where $Q$ is the quotient field of $R$. Since $Q$ is flat, the class $\mathcal {MC}$ of Matlis cotorsion modules contains the class $\mathcal {EC}:=\mathcal Fl^{\perp_1}$ 
of \emph{Enochs cotorsion modules}. Denoted by $\mathcal{TF}$ the class of torsion-free modules, the latter class $\mathcal {EC}$ contains the class $\mathcal {WC}:=\mathcal{TF}^{\perp}$ of Warfield cotorsion modules. Thus we have the following chain of cotorsion pairs, ordered with respect to the inclusion on the first class:
\[({}^{\perp_1}\mathcal {MC}, \mathcal {MC})\leq
(\mathcal Fl={}^{\perp_1}\mathcal {EC}, \mathcal {EC})\leq
(\mathcal{TF}={}^{\perp_1}\mathcal {WC}, \mathcal {WC}).\]
The modules in ${}^{\perp_1}\mathcal {MC}$ are called  \emph{strongly flat}.
The Enochs and Warfield cotorsion pairs $(\mathcal Fl, \mathcal {EC})$ and $(\mathcal{TF}, \mathcal {WC})$ are hereditary and the classes of flat and torsion free modules, as well known, define a homological dimension. In general the Matlis cotorsion pair $({}^{\perp_1}\mathcal {MC}, \mathcal {MC})$ is not hereditary and therefore strongly flat modules do not define a homological dimension;  precisely, the Matlis cotorsion pair is hereditary, and so strongly flat modules define a homological dimension, if and only if the quotient field $Q$ of $R$ has projective dimension $\leq 1$,  i.e. $R$ is a Matlis domain \cite[Section~10]{M}.
\end{EX}

\section{Generalizing the Gorenstein dimension}

Auslander in \cite{A} introduced  the notion of Gorenstein dimension for finite modules over a commutative noetherian ring. More precisely, let $R$ be a commutative noetherian ring; following \cite[Definition~1.1.2]{C}   we say that a finite  $R$-module $M$ belongs to the \emph{G-class} $G(R)$  if :
\begin{enumerate}
\item $\Ext{m}{R}MR=0$ for  $m>0$
\item $\Ext{m}{R}{\Hom_R(M,R)}R=0$ for $m>0$
\item the canonical morphism $\delta_M \colon M\to \Hom_R(\Hom_R(M,R), R)$, $\delta_M(x)(\psi)=\psi(x)$,  is an isomorphism.
\end{enumerate}  
Any finite  module admitting a  $G(R)$-resolution of length $n$ is said to have \emph{Gorenstein dimension} at most $n$.  
In \cite[Theorem 1.2.7]{C} it is shown that $G(R)$ defines  a homological dimension on the category of finite $R$-modules.

\vskip.3truecm
Given an abelian category $\mathcal A$, we denote by $\mathcal K(\mathcal A)$ (resp. $\mathcal K^+(\mathcal A)$, $\mathcal K^-(\mathcal A)$, $\mathcal K^b(\mathcal A)$) the homotopy category of unbounded (resp. bounded below, bounded above, bounded) complexes of objects of $\A$ and by $\D(\mathcal A)$ (resp. $\D^+(\mathcal A)$, $\D^-(\mathcal A)$, $\D^b(\mathcal A)$) the associated derived category. 
In the sequel with $\D^*(\mathcal A)$ or $\D^{\dag}(\mathcal A)$ we will denote any of these derived categories.

Consider a right adjoint pair of contravariant functors $(\Phi, \Psi)$ between the abelian categories $\A$ and $\B$, with the natural  morphisms $\eta$ and $\xi$ as unities. 
Following \cite[Theorem~5.1]{H}, to guarantee the existence of the derived functors $\R^*\Phi:\D^*(\mathcal A)\to \D(\mathcal B)$ and $\R^{\dag}\Psi:\D^{\dag}(\mathcal B)\to \D(\mathcal A)$, we assume the existence of triangulated subcategories $\mathcal P$ of $\mathcal K^*(\mathcal A)$ and $\mathcal Q$ of  $\mathcal K^{\dag}(\mathcal B)$ such that:
\begin{itemize}
\item every object of $\mathcal K^*(\mathcal A)$ and every object of $\mathcal K^{\dag}(\mathcal B)$ admits a quasi-isomorphism into objects of $\mathcal P$ and $\mathcal Q$, respectively;
\item if $P$ and $Q$ are exact complexes in $\mathcal P$ and $\mathcal Q$, then also $\Phi(P)$ and $\Psi(Q)$ are exact.
\end{itemize}

Given complexes $X\in\D^*(\mathcal A)$ and 
$Y\in\D^{\dag}(\mathcal B)$, we have
$\R^*\Phi X= \Phi P$ and $\R^{\dag}\Psi Y= \Psi Q$, where $P$ is a complex in $\mathcal P$ quasi-isomorphic to $X$, and $Q$ is a complex in $\mathcal Q$ quasi-isomorphic to $Y$.

The functor $\Phi$ has \emph{cohomological dimension $\leq n$} if, for each $A$ in $\mathcal A$, we have $H^i(\R^*\Phi A)=0$  for $|i|>n$.


An object $A$ in $\mathcal A$ is called \emph{$\Phi$-acyclic} if $H^i(\R^*\Phi A)=0$ for any $i\not=0$. Similarly, $\Psi$-acyclic objects in $\B$ are defined.
\begin{DEF}
We say that an object $A\in \mathcal A$ belongs to the class $\mathcal G_{{\Phi}{\Psi}}$ if   
\begin{enumerate}
\item $A$ is $\Phi$-acyclic;
\item $\Phi(A)$ is $\Psi$-acyclic
\item the morphism $\eta_A\colon A\to \Psi\Phi(A)$  is an isomorphism.
\end{enumerate}  
\end{DEF}

Note that, since the category of modules over a ring $R$ has enough projectives, the total derived functor $\R\Hom(-, R)$ always exists (see \cite{S}). Thus the class $\mathcal G_{{\Phi}{\Psi}}$ for the adjoint pair $(\Phi, \Psi)=(\Hom(-, R), \Hom(-, R))$ in the category of finite $R$-modules, coincides with the $G(R)$-class introduced above if $R$ is a commutative noetherian ring.

We want to prove that the class $\mathcal G_{{\Phi}{\Psi}}$ associated to the right adjoint pair $(\Phi, \Psi)$ always defines a homological dimension.

First  we   prove that  the $\mathcal G_{{\Phi}{\Psi}}$-dimension   can be computed using  the cohomology groups $H^i(\R^*\Phi)$. As a consequence it follows that,  when the category $\A$ has enough projectives, the $\mathcal G_{{\Phi}{\Psi}}$-dimension can be  compared  with the projective dimension  (conf. \cite[Proposition~1.2.10]{C}).

\begin{PROP}\label{prop:sup}
Let $A$ be an object in $\mathcal A$ of finite $\mathcal G_{{\Phi}{\Psi}}$-dimension. 
Then 
\begin{itemize}
\item[(a)] $\mathcal{G}_{{\Phi}{\Psi}}$-dim $A=\sup\{i: H^i(\R^*\Phi A)\}\not =0$
\item[(b)] If $\A$ has enough projectives, then $\mathcal{G}_{{\Phi}{\Psi}}$-dim $A\leq \pd A$
\end{itemize}
\end{PROP}
\begin{proof}
(a) Let $\mathcal G_{{\Phi}{\Psi}}$-dim $A=n$. Therefore there exists an exact sequence
$0\to G_{n}\to G_{n-1}...\to G_0\to A\to 0$ with $G_i\in\mathcal G_{{\Phi}{\Psi}}$, $i=0,1,..., n$. 
By shift dimension we get $H^i(\R^*\Phi A)=0$ for each $i>n$. If $\sup\{i:H^i(\R^*\Phi A)\not =0\}<n$, let $K$ be the cokernel of $G_{n}\to G_{n-1}$. We will prove that $K$ belongs to $\mathcal G_{{\Phi}{\Psi}}$ contradicting the assumption $\mathcal G_{{\Phi}{\Psi}}$-dim $A=n$. Indeed, $K$ is $\Phi$-acyclic
since $H^i(\R^*\Phi K)\cong H^{(i+n-1)}(\R^*\Phi A)=0$ for each $i>0$; 
applying $\Psi$ to the short exact sequence 
$0\to \Phi K\to \Phi G_{n-1}\to \Phi G_{n}\to 0$ and comparing it with the short exact sequence $0\to G_{n}\to G_{n-1}\to K\to 0$, we get that $\Phi K$ is $\Psi$-acyclic and the unity $\eta_K$ is an isomorphism. 

(b) If $\A$ has enough projectives, then any object $A$ in $\mathcal A$ admits a projective resolution $P$. Since the projectives are $\Phi$-acyclic, we have $\R^*\Phi  A=\Phi P$ and then \[\sup\{i:H^i(\R^*\Phi  A)\not =0\}=\sup\{i:H^i(\Phi P)\not =0\}\leq \pd A. \qedhere \]
\end{proof} 

Observe that, differently from the $G(R)$-dimension,  the inequality between the $\mathcal G_{{\Phi}{\Psi}}$-dimension and the projective dimension can be strict also  for objects of finite projective dimension (cf.  \cite[Proposition~1.2.10]{C}).

\begin{EX}
Let $\Lambda$ be the path algebra of  the quiver 
$$
\xymatrix{ 1 \ar[r]
 & 2  \ar[r]
 & 3}
$$
 Let us consider the  module $_{\Lambda}U=\begin{smallmatrix} 1\\2\\3 \end{smallmatrix}\oplus \begin{smallmatrix} 2\\3 \end{smallmatrix}\oplus \begin{smallmatrix} 2 \end{smallmatrix}$ and let $S=\End_{\Lambda}(U)$.
  Consider the adjoint pair $(\Hom_{\Lambda}(-, U),  \Hom_{S}(-, U))$: since $\Ext{1}{\Lambda}UU=0$, $\Ext{1}{S}SU=0$ and $U\cong \Hom_{S}(\Hom_{\Lambda}(U, U), U)$, the $\Lambda$-module $U$ belongs to  $\mathcal G_{{\Phi}{\Psi}}$, where $(\Phi, \Psi)=(\Hom_{\Lambda}(-, U), \Hom_S(-, U))$. Thus $U$    has projective dimension one, but obviously $\mathcal G_{{\Phi}{\Psi}}$-dimension 0.
\end{EX}

 In order  to prove that the class  $\mathcal G_{{\Phi}{\Psi}}$ defines a homological dimension, we also need to recall some notions and results on derived categories.   By  \cite[Lemma~13.6]{K}  we know that, in our assumptions, $(\R^* \Phi, \R^{\dag} \Psi)$ is a right adjoint pair in the derived categories $\D^*(\A)$ and $\D^{\dag}(\B)$, with unities $\hat\eta$ and $\hat\xi$ naturally inherited from the unities $\eta$ and $\xi$. In \cite{MT} a complex $X\in \D^*(\A)$ is called \emph{$\D$-reflexive} if the morphism $\hat\eta_X$  is an isomorphism in  $\D^*(\A)$. An object $A\in\A$ is called \emph{$\D$-reflexive} if it is $\D$-reflexive  as a stalk complex.
 
 
 
 \begin{LEMMA}\label{lemma:precedente}
 Let $X\in \A$ such that $X$ is $\Phi$-acyclic and $\Phi(X)$ is $\Psi$-acyclic. Then $\hat\eta_X$ is a quasi-isomorphism if and only if $\eta_X$ is an isomorphism. In particular any object in $\mathcal G_{{\Phi}{\Psi}}$ is $\D$-reflexive.
 \end{LEMMA}
 \begin{proof}  In general,  if $C\in \D^{*}(\A)$ and $L$ is a complex quasi-isomorphic to $C$ such that any term $L_i$   of $L$ is   $\Phi$-acyclic and $\Phi(L_i)$ is $\Psi$-acyclic, then $\hat{\eta}_C$ coincides with  $\eta_L$, where $\eta_L$ is the term-to-term extension of the unity $\eta$ to the triangulated category  $\mathcal{K}^{*}(\A)$ (cf. \cite{MT}). Then we easily get the statement.
 \end{proof}

 \begin{COR}\label{cor:Dreflexive}
 Any object $A$ in $\mathcal A$ of finite $\mathcal G_{{\Phi}{\Psi}}$-dimension is $\D$-reflexive.
 \end{COR}
 \begin{proof}
 Let $\mathcal G_{{\Phi}{\Psi}}$-dim $A=n$. Therefore there exists an exact sequence
$0\to G_{n}\to G_{n-1}...\to G_0\to A\to 0$ with $G_i\in\mathcal G_{{\Phi}{\Psi}}$, $i=0,1,..., n$.  Therefore in the bounded derived category $\D^b(\mathcal A)$, $A$ is quasi-isomorphic to the complex $G:=0\to G_{n}\to G_{n-1}...\to G_0\to 0$. Since $G$ is a complex with $\D$-reflexive terms by Lemma~\ref{lemma:precedente}, we conclude by \cite[Theorem~3.1,(1)]{MT} that $A$ is $\D$-reflexive.
 \end{proof}
  
 \begin{PROP}\label{prop:?}
 If $X\in\A$ is $\Phi$-acyclic and $\D$-reflexive, then $X$ belongs to $\mathcal G_{{\Phi}{\Psi}}$.
  \end{PROP}
\begin{proof}  Since $X$ is $\Phi$-acyclic, $\R^*\Phi X$ is quasi isomorphic to the stalk complex $\Phi(X)$. Moreover, for $X$ is $\D$-reflexive, we get that $\R^{\dag} \Psi  (\Phi X)\cong \R^{\dag} \Psi  (\R^*\Phi X) $ is quasi-isomorphic to $X$. Thus  $\Hm^i(\R^{\dag} \Psi  (\Phi X))=0$ for any $i\neq 0$ and so $\Phi X$ is $\Psi$-acyclic. Finally we conclude since, by the previous lemma,  $\eta_{X}$  is an isomorphism. 
\end{proof}

\begin{TEO}\label{thm:homdim} 
The class $\mathcal G_{{\Phi}{\Psi}}$ defines a homological dimension.
\end{TEO}
\begin{proof}
 Let us consider a long exact sequence $0\to G_{n} \to G_{n-1} \to \dots \to G_0 \to X\to 0$ with $G_i\in \mathcal G_{{\Phi}{\Psi}}$.
By Corollary~\ref{cor:Dreflexive}, $X$ is $\D$-reflexive.
 Consider  now a long exact sequence $0\to X_{n} \to F_{n-1} \to \dots \to F_0 \to X\to 0$ with $F_i\in \mathcal G_{{\Phi}{\Psi}}$.
 The $\D$-reflexive objects are a thick subcategory of $\A$  (see \cite{MT}), i.e, if two terms of a short exact sequence in $\mathcal A$ are $\D$-reflexive,  then also the third is $\D$-reflexive. Therefore, by induction,  it follows that $X_{n}$ is $\D$-reflexive. Since by Proposition~\ref{prop:sup} $H^i(\R^*\Phi X)=0$ for each $i>n$, by shift dimension $X_{n}$ is $\Phi$-acyclic, and so we conclude that $X_{n}$ belongs to $\mathcal G_{{\Phi}{\Psi}}$.
\end{proof}

In \cite{MT}, the authors were interested in characterizing the $\D$-reflexive objects associated to  a given adjoint pair $(\Phi, \Psi)$. Assume $\mathcal A$ is a module category and denote by $FP_n$ the class of modules $A$ which have  an exact resolution
\[P_{n-1}\to ...\to P_1\to P_0\to A\to 0,\]
where the $P_i$'s are finitely generated projectives. In particular $FP_1$ is the class of finitely generated modules. Then the $\D$-reflexive modules in $FP_n$ can be characterized through their $\mathcal G_{{\Phi}{\Psi}}$-dimension.
%

\begin{TEO}\label{thm:fingen} Let $\mathcal A=R\lmod$ for an arbitrary ring $R$. Assume $_{R}R$ to be $\D$-reflexive and  $\Phi$ of cohomological dimension $\leq n$ .
Then   a module $M\in FP_n$ is $\D$-reflexive if and only if it has $\mathcal G_{{\Phi}{\Psi}}$-dimension $\leq n$.
\end{TEO}
\begin{proof}  The sufficiency of the finiteness of the $\mathcal G_{{\Phi}{\Psi}}$-dimension is proved in Corollary~\ref{cor:Dreflexive}.
Conversely,  suppose $M$ to be a $\D$-reflexive module in $FP_n$. 
Let $0\to K\to P_{n-1}\to \dots \to P_0\to M\to 0$ be an exact sequence with the $P_i$'s finitely generated projectives.  Since $_{R}R$ is assumed to be $\D$-reflexive, any $P_i$ is $\D$-reflexive, and so we get that $K$ is $\D$-reflexive. Since $\Phi$ has cohomological dimension $\leq n$,  by shift dimension we get that $K$ is $\Phi$-acyclic. Then, by Proposition~\ref{prop:?}, we conclude that  $K$ belongs to $\mathcal G_{{\Phi}{\Psi}}$.
\end{proof}

\end{document}